\documentclass{amsart}[12pt]
\usepackage{amscd}

\usepackage{amsmath}
\usepackage{amssymb}
\newcommand{\R}{{\mathbb R}}
\newcommand{\C}{{\mathbb C}}
\newcommand{\Z}{{\mathbb Z}}

\newcommand{\tJ}{\tilde{J}}

\newcommand{\ah}{\mathcal{A}_H}
\newcommand{\rs}{{\mathbb R} \times S^1}
\newcommand{\hk}{H_{\lambda} \natural K_{\delta}}
\newcommand{\hl}{\bar{H}_{\lambda}}
\newtheorem{theo}{Theorem}
\newtheorem{lm}{Lemma}
\newtheorem{df}{Definition}
\newtheorem{prop}{Proposition}

\newtheorem{rema}{Remark}
\newtheorem{prob}{Problem}
\begin{document}

\title[ Symplectic rigidity, fixed points and global perturbations.]{Symplectic rigidity, symplectic fixed points and global
perturbations of Hamiltonian systems.}

\author{Dragomir L. Dragnev \\}

\address{Courant Institute of Mathematical Sciences\\
New York University \\ 251 Mercer street \\New York, NY 10012}
\email{dragnev@cims.nyu.edu}
\thanks{Preliminary version. November, 2005}
\begin{abstract}
In this paper we study a generalized symplectic fixed point
problem, first considered by J. Moser in \cite{M}, from the point
of view of some relatively recently discovered symplectic rigidity
phenomena. This problem has interesting applications concerning
global perturbations of Hamiltonian systems.
\end{abstract}

\maketitle
%%%%%%%%%%%%%%%%%%%%%%%%%%%%%%%%%%%%%%%%%%%%%%%%%%%%%%%%%%%%%%%%%%%%%%%%%%%
\section{Introduction} Let $(M, \omega)$ be a symplectic manifold.
We recall that a submanifold $N$ of $M$ is called
\textit{coisotropic} if at any point $x \in N$ we have that $(T_x
N)^{\omega} \subseteq T_x N$, where $(T_x N)^{\omega} = \{v \in
T_x M \mid \omega(v,u) = 0, u \in T_x N\}$. The distribution,
$(TN)^{\omega}$, on $N$ is integrable (see \cite{MS,M}), and
therefore gives rise to a foliation on $N$. Denote by
$\mathcal{L}_x N$ the leaf of this foliation through $x \in N$. We
are interested in the following geometric problem.
\begin{prob} Given a symplectomorphism $\phi$ of $M$ (i.e., $\phi ^{*}
\omega =\omega$), under what conditions on $\phi$ and possibly on
$N$, there exists a point $x \in N$ so that its image $\phi(x)$
lies on a leaf through $x$, i.e., $\phi(x) \in \mathcal{L}_x N$.
\end{prob}

In this paper we are going to study the above problem for a
special class of coisotropic submanifolds, following Ph. Bolle,
\cite{B}, we present the following:
\begin{df} \label{cont} Let $N$ be a $k$-codimensional compact coisotropic
submanifold of a symplectic manifold $(M^{2n},\omega)$ and $1 \leq
k\leq n$. $N$ is called ``$k$-contact" if there exist $k$ 1-forms
$\alpha_1, \ldots, \alpha_k$ defined on $N$ so that
\begin{enumerate}
\item $d \alpha_i = \omega_{|N}$ for $i=1,\ldots, k$. \item For
all $x \in N$, $\alpha_1 \wedge \ldots \wedge \alpha_k \wedge
\omega ^{n-k} (x) \neq 0$.
\end{enumerate}
Equivalently one can state the second condition as follows, for
all $x \in N$, the restrictions of $\alpha_1 (x), \ldots, \alpha_k
(x)$ to $Ker \omega_{|N}$ are linearly independent.
\end{df}
Our main result is:
\begin{theo} \label{D} Let $N$ be a compact submanifold of $(\R^{2n},
\omega_0)$ of $k$-contact type. Let $\phi$ be the time-1 map of a
compactly supported Hamiltonian $H$ on $[0,1] \times \R^{2n}$ such
that $E(\phi) < c_{FH} (N)$. Then there exists $x \in N$ such that
$\phi(x) \in \mathcal{L}_x N$.
\end{theo}
Here $c_{FH}$ stands for the Floer-Hofer capacity as defined in
\cite{He1}, see Section \ref{FlHo}, $\omega_0 = -d\lambda_0$ with
$\lambda_0 = 1/2 \sum_{j=1}^n (y_j dx_j -x_jdy_j)$, is the
standard symplectic structure on $\R^{2n}= \C^n$, and the energy
$E(\phi)$ is defined as follows. Denote by $\mathcal{F}$ the space
of all smooth functions $H: [0,1] \times M \longrightarrow \R$
with compact support. To every such function one can associate a
symplectic map $\phi_{H}= \varphi^1$, where $\varphi^{t}$ is the
flow of the Hamiltonian vector field, $X_{H_{t}}$, defined by the
equation $\omega(H_t, \cdot) = -dH(\cdot)$. We call a symplectic
map $\phi$ - Hamiltonian if $\phi = \phi_{H}$ for some function $H
\in \mathcal{F}$. Following Hofer, \cite{H}, we define the norm of
$H$ to be $\|H\| = \sup H - \inf H$ and the energy of a
Hamiltonian map,
\begin{equation}E(\phi) = \inf _{H\in \mathcal{F}}\{ \|H\| \mid \phi =
\phi_{H}\}
\end{equation}
Problem 1 was first considered in \cite{M} and J. Moser proved
that Problem 1 has a solution if $(M, \omega = d \alpha)$ is a
simply connected, exact symplectic manifold, $N$ is a compact
coisotropic submanifold of $M$ and $\phi$ is an exact
symplectomorphism of $M$, (that is $\phi ^* \alpha - \alpha$ is
exact), which is $\mathrm{C}^1$ close to the identity. Obviously
Moser's result is of local nature. In 1989, I. Ekeland and H.
Hofer derived more global versions of this theorem, for the case
where $N$ is a compact hypersurface of restricted contact type in
$(\R^{2n}, \omega_0)$ and the map $\phi$ is a Hamiltonian
symplectomorphism. They presented various conditions on the map
$\phi$ for which the problem above has a solution, see \cite{EH1}
for details. Here we recall that a hypersurface $N$ in a
symplectic manifold $(M^{2n}, \omega)$ is called of contact type
if it is $1$-contact in terms of Definition \ref{cont}. $N$ is
said to be of restricted contact type if, in addition, the form
$\alpha_1$ can be extended to $M$ satisfying $d\alpha_1 = \omega$.
In \cite{H}, H. Hofer, proved a very surprising result stating
that the problem above has a solution if $N$ is a compact
hypersurface of restricted contact type in $(\R^{2n}, \omega_0)$,
and $\phi$ is a time-1 map of a compactly supported Hamiltonian,
provided the energy of $\phi$ is bounded by the Ekeland-Hofer
capacity of $N$, (defined in \cite{EH2, EH3}), i.e. $$E(\phi) \leq
c_{\textrm{EH}}(N).$$ Theorem \ref{D} extends, in a way, Hofer's
result to coisotropic submanifolds of higher codimension and even
in codimension one we do not assume $N$ to be of restricted
contact type. We point out that the results in \cite{EH1, H} were
obtained by using variational methods which are somehow restricted
to the Euclidian case. Another limitation of the variational
approach, even in the Euclidean case, is that it does not allow us
to gain the needed control of the gradient trajectories of the
Hamiltonian action functional, defined in Section \ref{symh},
(\ref{eq10}). On the other hand, mixing the variational approach
with pseudo-holomorphic curve methods in the spirit of Floer
homology, allows us to regain this control from a geometric or
rather topological prospective. Namely the idea behind the proof
of the main theorem is to foliate a small neighborhood of $N$ into
diffeomorphic images of $N$. Then we consider the critical points
of a special action functional and establish the existence of a
critical point which is a closed trajectory for a special
Hamiltonian and consists of two arcs one is $\psi^t x$ s.t.
$\psi^1 = \phi$ and the other arc connects $x$ and $\phi (x)$
through a path which is on the leaf through $x$ on a nearby image
of $N$. We do this by studying the symplectic homology groups of
this neighborhood. The existence of the closed trajectory of the
type described above is a consequence of the non-vanishing of
certain Floer homology groups filtered by the action. Taking
smaller and smaller neighborhoods of $N$, and repeating the
previous step we get a family of closed trajectories of this type
and we want to take a limit of these which will be a solution of
Problem 1. The subtle part is to show that the lengths of the arcs
which are on the leaves of the nearby images of $N$ are uniformly
bounded. We achieve this by getting some bounds on the action of
the critical point which comes automatically from the fact that we
work with filtered Floer homology groups plus some additional
information coming from the functorial properties of the
symplectic (Floer) homology i.e. that the critical point is a
deformation of the constant solution of a certain Hamiltonian.
This information is impossible to be detected by the variational
approach and that is the reason, which directs us to work with the
Floer-Hofer capacity which is based on the symplectic homology of
Floer and Hofer, \cite{FH1}. Perhaps it is worth mentioning that
the contact condition is significant for the Hofer's theorem. We
refer to a recent paper of V. Ginzburg, \cite{G}, for a discussion
about the significance of the contact condition in the various
existence and almost existence results of periodic orbits on
hypersurfaces and its importance for the validity of the Weinstein
conjecture. Since the methods we are going to employ are
reminiscent to some of the methods used to prove certain cases of
the Weinstein conjecture, we must impose some sort of a contact
type condition on $N$, and this justifies our choice of the
$k$-contact condition. We postpone the discussion on what are the
right conditions on $N$ and the consideration of Problem 1 on more
general symplectic manifolds, most notably cotangent bundles, to
our forthcoming paper \cite{D}.

As an almost immediate application of Theorem \ref{D} we consider
the Hamiltonian system describing the motion of $n$ independent
harmonic oscillators on $\R^{2n}$, with Hamiltonian
\begin{equation} \label{harmosc}
H_0= \frac{1}{2} \sum_{j=1}^n m_j (x_j ^2 + y_j ^2).
\end{equation}
It is well-known that this system is integrable with first
integrals $G_j = x_j^2 + y_j ^2$. Consider for suitable positive
constants $c, c_1, \ldots , c_{k-1}$, where $2\leq k \leq n$, the
level manifold
$$N_{(c, c_1, \ldots , c_{k-1})} = \{H_0=c, G_j = c_j, j = 1, \ldots,
k-1\}.$$ It is not hard to see that $N_{(c, c_1, \ldots ,
c_{k-1})}$ is a compact, coisotropic, $k$-codimensional
submanifold of $(\R^{2n}, \omega_0)$, see \cite{M}. In fact we
shall see that it is of $k$-contact type. In polar coordinates
$$ x_j -i y_j = r_j e^{i \theta_j}$$ one has $$r_j^2 = c_j$$
\begin{equation}\label{cond}
\frac{1}{2} \sum_{j=k}^{n} m_j r_j^2 = c - \frac{1}{2}
\sum_{j=1}^{k-1} m_j c_j > 0.
\end{equation}
The flow generated by $G_j$ is given by $ r_j \to r_j$ and
$\theta_j \to \theta_j - \delta_{jl} \tau_l$, where $\delta_{jl}$
is the Kronecker symbol. The leaves through a point
$(r^*,\theta^*)$ are given by $$ r_j = r_j^*; \theta_j =
\theta_j^* + \sum_{l=1}^{k-1} \delta_{jl} \tau_l +m_j \tau_{k}$$
where $j=1, \ldots , n$ and $\tau_1, \ldots, \tau_k$ are the $k$
parameters on the leaf.

Now let us consider a nonautonomous, compactly supported,
Hamiltonian perturbation $H_1(t,x): \R \times \R^{2n} \to \R$,
such that $supp H_1 \subset [0,1] \times K$, for some compact
subset $K$ of $\R^{2n}$. We have the following theorem.

\begin{theo}\label{DD} Assume that $$\parallel H_1 \parallel < \min \bigg\{
\min_{p=1,\ldots, k-1} \{\pi c_p \}, \min_{p=k,\ldots, n} \{\pi
\frac{(2c-\sum_{j=1}^{k-1} m_j c_j)}{m_p}\}\bigg\}$$ then there
exists a solution $y$ of the perturbed system with Hamiltonian
$H_0 + H_1$ which aside from phase shifts $\tau_1,\ldots, \tau_k$
of $\theta_1, \ldots, \theta_{k-1}, t$ returns to the continuation
of the unperturbed orbit. In particular the integrals $H_0,
G_1,\ldots, G_{k-1}$ have the same value for $t \in
(-\infty,0)\bigcup (1,\infty)$.
\end{theo}

The paper is organized as follows. In Section \ref{symh} we review
the definition and some of the properties of the Floer homology
and the symplectic homology respectively as well as the definition
of the symplectic capacities and the Floer-Hofer capacity in
particular. In Section \ref{cd} we review some consequences of
Definition \ref{cont}. The proof of the main theorem is done in
Section \ref{mt} and Section \ref{appl} is devoted to the proof of
Theorem \ref{DD}.

%%%%%%%%%%%%%%%%%%%%%%%%%%%%%%%%%%%%%%%%%%%%%%%%%%%%%%%%%%%%%%%%%%%%%%%%%%%%%%%%%%%%%%%%%%

\section{Consequences of the contact definition.} \label{cd} In this section
we review some useful results from \cite{B}. Denote by
$B_{\varepsilon}^k$ the ball with center $0$ and radius
$\varepsilon$ in $\R^k$. Then we have the following lemma.

\begin{lm} \label{lm1} Let $N$ be a smooth, compact, connected coisotropic
submanifold of a symplectic manifold $(M, \omega)$ which is of
$k$-contact type. Then there exists $\varepsilon >0$, an open
neighborhood $U$ of $N$ in $M$ and a diffeomorphism $\psi : N
\times B_{\varepsilon}^k  \rightarrow U$ such that:

i) For all $x \in N$ we have $\psi(x, 0) =x$;

ii) $\psi^* \omega = (1 + \sum_{j=1} ^k y_j) q^*(\omega_{|N}) +
\sum_{j=1} ^k dy_j \wedge q^*(\alpha_j)$;\\ where the $1$-forms
$\alpha_j$ are the ones from Definition \ref{cont}, $q : N \times
B_{\varepsilon}^k  \rightarrow N$ is the projection onto the first
factor and $y_1, \ldots, y_k$ are coordinates in
$B_{\varepsilon}^k$.
\end{lm}
We have some useful consequences. We set the notation,
$$ r = q \circ \psi ^{-1} : U \rightarrow N$$
\begin{equation} \label{beta} \beta _j = \psi ^{-1 *} (q ^* (\alpha_j)) = r^*
\alpha_j
\end{equation}
$$W = r^*(\omega_{|N})$$
\begin{equation} z_j = y_j \circ \psi^{-1}
\end{equation}

With this notation we have from the lemma that in $U$ the
following is true:
\begin{equation}\omega = (1 + \sum_{j=1} ^k z_j) W +
\sum_{j=1} ^k dz_j \wedge \beta_j
\end{equation}
Denote by $X_{z_j}$ the Hamiltonian vector field associated to
$z_j$, i.e., $\omega(X_{z_j},\cdot) = -dz_j(\cdot)$. Then by the
above lemma we have
\begin{equation}\label{eq1} dz_j(X_{z_i}) = 0
\end{equation}
\begin{equation}\label{bet} \beta_j(X_{z_i}) = \delta_{ij}
\end{equation}
\begin{equation} W(X_{z_i}, \cdot)=0
\end{equation}
for $i,j = 1, \ldots k$. It follows from (\ref{eq1}), that the
functions $z_1, \ldots, z_k$ are in involution. Set for $\nu =
(\nu_1, \ldots, \nu_k) \in B_{\varepsilon}^k$, $N(\nu) =
\cap_{j=1}^k z_j ^{-1} (\nu_j)$. Then $r_{|N(\nu)}$ is a
diffeomorphism from $N(\nu)$ onto $N$. Moreover it is not hard to
see that $N(\nu)$ is a coisotropic submanifold of $M$ and
$(TN(\nu))^{\omega}$ is spanned by $X_{z_1}, \ldots, X_{z_k}$.
From this follows that if we have a trajectory $x(t)$ satisfying
the equation
\begin{equation} \label{leaf} \dot{x} = \sum_{j=1} ^k \gamma_j X_{z_j} (x(t))
\end{equation}
for some coefficients $\gamma_j$, then $x(t)$ will be on the leaf
through $x(0)$ of $N(\nu)$ where $\nu=(z_1(x(t)), \ldots ,
z_k(x(t))$. This observation will play a significant role in the
proof of Theorem \ref{D}. We conclude this section by noticing
(due to Lemma \ref{lm1}), that we can foliate a neighborhood of a
$k$-contact submanifold $N$ into coisotropic images of $N$ in $M$.

%%%%%%%%%%%%%%%%%%%%%%%%%%%%%%%%%%%%%%%%%%%%%%%%%%%%%%%%%%%%%%%%%%%%%%%%%%%%%%%%%%%%%%%%%%

\section{Review of the symplectic homology and the definition of $c_{FH}$.}
\label{symh} In this section we review briefly the definition and
the properties of the Floer homology and the symplectic homology.
Based on the properties of the symplectic homology we will present
a very useful symplectic invariant called the Floer-Hofer
capacity, in the terminology of D. Hermann, \cite{He1}.

%%%%%%%%%%%%%%%%%%%%%%%%%%%%%%%%%%%%%%%%%%%%%%%%%%%%%%%%%%%%%%%%%%%%%%%%%%%%%%%%%%%%%%%

\subsection{Floer Homology for the Hamiltonian action functional.}
The Floer homology is an infinite-dimensional equivalent to the
Morse theory. In other words it can be thought as a version of
Morse theory for the Hamiltonian action functional. Here we recall
the definition and the properties of Floer homology. Details can
be found in \cite{HZ,Sa} or in the A. Floer original paper
\cite{F}.

Let $(M^{2n}, \omega)$ be a closed symplectic manifold, which is
symplectically aspherical, that is,
$$\omega |_{\pi_2(M)} =0 \textrm{  and  } c_1(TM) |_{\pi_2(M)}
=0$$ where $c_1(TM)$ is the first Chern class of the tangent
bundle of $M$. Let $H \in \mathcal{F}$ be a time-dependent
function on $M$ and $X_{H}$ be its Hamiltonian vector field.
Denote by $\mathcal{P}(H)$ the set of contractible one-periodic
orbits of $X_H$. Let $\Lambda M$ be the space of smooth
contractible loops in $M$. We define the Hamiltonian action
functional, $\mathcal{A}_H : \Lambda M \to \R$, associated with $H
\in \mathcal{F}$, as follows,
\begin{equation}\label{eq10}
\mathcal{A}_H (x) = \int_{D} u^{*} \omega - \int_0 ^1 H(t,x(t)) dt
\end{equation}
with $D$ being the closed unit disc ($\partial D = S^1$), and $u :
D \to M$ an extension of $x$ so that $u|_{\partial D} = x$. This
functional is well defined because of our assumption that $M$ is
symplectically aspherical. As we mentioned above the Floer
homology may be viewed as a Morse theory on $\Lambda M$. To be
precise we denote by $\mathcal{J}_M$ the space of
$\omega$-compatible almost complex structures on $M$, i.e. the
space of all $J: TM \to TM$ such that $J^2 = - Id$ and
\begin{equation}\label{eq11}
\omega (\xi, J \eta) = g_J (\xi, \eta) \textrm{  for all  } \xi,
\eta \in TM
\end{equation}
so that $g_J$ is a Riemannian metric on $M$. Now consider the
$L_2$-metric induced on $\Lambda M$ by $g_J$. Then the gradient
of $\mathcal{A}_H$ is given by
\begin{equation}\label{eq12}
\nabla_J \mathcal{A}_H (x) = -J \dot{x} - \nabla H(t,x)
\end{equation}
In view of (\ref{eq12}) and the fact that $X_H = J \nabla H$, we
notice that the critical points of $\mathcal{A}_H$, are exactly
the one-periodic solutions of the Hamiltonian equations $\dot{x}=
X_H (x), x(0)=x(1)$ i.e. the elements of $\mathcal{P}(H)$. The set
of critical values of $\mathcal{A}_H$ is called the action
spectrum of $H$ and denoted by $\Sigma (H)$. Of course,
(\ref{eq12}) does not define a flow on $\Lambda M$ but despite
that we are going to consider the gradient lines of $\nabla
\mathcal{A}_H$ as the solutions of the following elliptic equation
of Cauchy-Riemann type:
\begin{equation}\label{eq13}
\frac{\partial u}{\partial s} + J(t,u) \frac{\partial u}{\partial
t} + \nabla H(t,u) = 0 \textrm{  for  } u \in C^{\infty}( \R
\times S^1, M)
\end{equation}
Given two critical points $x^+ , x^- \in \mathcal{P}(H)$ of $\ah$
we consider the space of solutions $\mathcal{M}(x^- ,x^+,J,H)$ of
(\ref{eq13}) connecting $x^-$ and $x^+$,
$$
\mathcal{M}(x^- ,x^+ ,J,H) = \{ u \in C^{\infty}( \R \times S^1,
M) | \textrm{ (\ref{eq13}) and } \lim_{s \to \pm \infty} u(s,t) =
x^{\pm} (t) \}$$ An element, $u$, of $\mathcal{M}(x^- ,x^+ ,J,H)$
will be called a \textit{Floer trajectory}. In this situation the
difference of the actions between the ends is given by the energy,
$E_J(u)$ of the Floer trajectory $u$, defined as follows,
\begin{equation}\label{eq14}
\ah ( x^+) - \ah (x^-) = \int_{\rs} g_J (\frac{\partial
u}{\partial s}, \frac{\partial u}{\partial s}) ds dt = E_J (u)\geq
0
\end{equation}
Notice that the action is increasing along the gradient
trajectory, that is, $\frac{\partial \ah(u(s,\cdot))}{\partial s}
= \| \nabla \ah (u( s, \cdot)) \|^2 \geq 0$. It is not hard to see
that if $E_J(u) = 0$, then $u$ is independent of $s$, one-periodic
solution of the hamiltonian equations for $H$. Assume that the
elements of $\mathcal{P}(H)$ are non-degenerate, which means that
if $x(t) = x(t+1) \in \mathcal{P}(H)$, then
$$\det (Id - d \varphi_{H} ^1 (x(0)) \neq 0$$
where $\varphi_{H} ^t$ is the flow of $X_H$. With this assumption
and utilizing our assumption that $c_1 (TM) |_{\pi_2 (M)} = 0$,
the elements of $\mathcal{P}(H)$ are graded by their
Conley-Zehnder index, $\mu _{CZ}$; see \cite{SZ}. The key result
concerning moduli spaces $\mathcal{M}(x^- ,x^+ ,J,H)$ is the
following, see \cite{Sa},
\begin{theo}
For generic choices of $J$ and $H$, the moduli spaces
$\mathcal{M}(x^- ,x^+ ,J,H)$ are compact, finite dimensional,
manifolds, of dimension $\mu _{CZ} (x^+) - \mu _{CZ} (x^-)$.
\end{theo}
Following Floer, we define the Morse-Witten complex associated
with $H$ as a graded $\Z_2$-vector space
$$CF (H) = \bigoplus _{x \in \mathcal{P}(H)} \Z_2 x$$
We proceed by defining the \textit{Floer boundary operator}, by
$$\partial ^{H,J} x = \sum_{ y \in \mathcal{P}(H); \mu_{CZ}(x) -
\mu_{CZ} (y) = 1} \nu (x,y) y $$ where $\nu (x,y)$ stands for the
number (mod 2) of the elements in $\bar{\mathcal{M}}(x^- ,x^+
,J,H) = \mathcal{M}(x^- ,x^+ ,J,H) / \R$. In the last expression,
we observe that $\R$ acts freely by translation on the Floer
trajectories, and we mod out its action. The operator $\partial
^{H,J}$, satisfies $\partial ^{H,J} \circ
\partial ^{H,J} = 0$, thus allowing us to define the Floer
homology groups,
$$HF_* (H,J) = Ker \partial ^{H,J} / Im \partial ^{H,J}$$
It turns out that these groups are independent of the generic
choice of $J$, $HF_* (H) = HF_* (H,J)$. Later on it will be useful
to consider the Floer homology groups filtered by the action and
we take a moment to review their construction. Let $-\infty < a
\leq b < \infty$ be two numbers, so that $a,b \notin \Sigma (H)$.
Then we define for $a$, (respectively $b$), $\mathcal{P}^a (H) =
\{ x \in \mathcal{P}(H) | \ah (x) < a \}$ and $$CF ^a (H) =
\bigoplus _{x \in \mathcal{P}^a (H)} \Z_2 x$$ Then $CF^a (H)$ is a
subcomplex of $CF^b (H)$ and we consider the quotient complex
$CF^{[a,b)} (H) = CF^a (H) / CF^b (H)$. The filtered Floer
homology groups, $HF^{[a,b)} (H)$, are the homology groups of
$CF^{[a,b)} (H)$ with the induced boundary operator.
%%%%%%%%%%%%%%%%%%%%%%%%%%%%%%%%%%%%%%%%%%%%%%%%%%%%%%%%%%%%%%%%%%%%%%%%%%%%%%%%%%%%%%%%%%

\subsection{Symplectic Homology.}\label{SH}There are several different versions of
the symplectic homology. Originally it was introduced by A. Floer
and H. Hofer for bounded, open sets in $\R ^{2n}$, \cite{FH1}, by
further developing the idea behind the Floer theory and combining
that with ideas of I. Ekeland and H. Hofer about using the
Hamiltonian dynamics to study the symplectic rigidity,
\cite{EH1,EH2}. Later on versions of the symplectic homology,
concerning relatively compact sets in symplectic manifolds with
contact type boundary, \cite{CFH} and symplectic manifolds with
contact type boundary, \cite{V}, were developed. Here we are going
to use the original version of the symplectic homology from
\cite{FH1}, with $\Z_2$-coefficients, and refer the interested
reader to the survey paper of A. Oancea, \cite{Oa}, where the
different versions of the symplectic homology are compared.

Let $U$ be a bounded open set in $(\R^{2n} = \C^n, \omega_0)$.
Next we define the set of \textit{admissible} Hamiltonian
functions, $\mathcal{H}_{ad}(U)$.
\begin{df} A function $H: S^1 \times \C^n \to \R$ is called
\textit{admissible}, $H \in \mathcal{H}_{ad}(U)$.
\begin{enumerate}
\item $H|_{\bar{U}} <0$ for all $t \in S^1$; \item There is a
positive-definite matrix $A$ so that $\frac{|H'(t,z) - Az|}{|z|}
\to 0$ as $|z| \to \infty$, uniformly for $t \in S^1$; \item there
is a constant $c>0$ so that $$\|H'' (t,z) \| \leq c$$
$$|\frac{\partial}{\partial t} H(t,z)| \leq c(1+|z|)$$
\item the system $-i \dot{v} = Av$, admits no nontrivial
1-periodic solutions.
\end{enumerate}
\end{df}
Before we proceed, let us comment on the conditions in the above
definition. The first condition restricts a function $H \in
\mathcal{H}_{ad}(U)$ on the set $U$, and $H$ is allowed to
increase fast near the boundary of $U$. The second condition,
determines the asymptotic behavior of $X_H$, which combined with
the fourth condition , allows us to conclude that all 1-periodic
orbits of $H$ are contained in a compact set together with their
connecting (Floer) trajectories. The third condition is for
technical purposes and allows one to do the necessary estimates
needed for the well-definedness of the Floer homology in this
situation, i.e. in the case of open symplectic manifolds. Denote
by $\mathcal{H}_{reg} (U)$ the set of admissible Hamiltonians with
non-degenerate 1-periodic orbits and by $\mathcal{J}$ the set of
almost complex structures, compatible with the standard symplectic
structure $\omega_0$, which are equal to the standard complex
structure $i$ outside of a compact set. In \cite{FH1}, the
transversality of the Floer's equation,(\ref{eq13}), is
established for a dense subset of $\mathcal{H}_{reg} (U) \times
\mathcal{J}$. Following the discussion of the previous section,
one can define the Floer homology groups, filtered by the action,
for a regular pair $(H,J)$. Symplectic homology arises from
certain functorial properties of Floer homology. Given regular
pairs $(H_1,J_1)$ and $(H_2, J_2)$, such that $H_1 \leq H_2$, on
$S^1 \times \C^n$, we consider a monotone homotopy connecting
them. That is a homotopy $(L(s,t,z), \tilde{J}(s,t,z))$ such that:
\begin{itemize}
\item $(L(s,t,z), \tilde{J}(s,t,z)) = (H_2(t,z), J_2(t,z))$ for
$s\leq -s_0$; \item $(L(s,t,z), \tilde{J}(s,t,z)) = (H_1(t,z),
J_1(t,z))$ for $s\geq s_0$; \item $\frac{\partial L}{\partial s}
\leq 0$ on $\R \times S^1 \times \C^n$; \item There is a smooth
path $A(s)$ of positive matrices so that $A(s) = A(-s_0)$ for $s
\leq -s_0$ and $A(s) = A(s_0)$ for $s \geq s_0$ and $$\lim _{|z|
\to \infty} \frac{|L'(s,t,z) - A(s)z|}{|z|} \to 0$$ plus we
require that if the system $-i \dot{v} = A(s)v$ has a non-trivial
1-periodic solution for some $s=s'$ then $\frac{d}{ds}|_{s=s'}
A(s)$ is positive definite.
\end{itemize}
Consider the parametrized version of the Floer equations
(\ref{eq13}),
\begin{equation}\label{monh12}
\frac{\partial u}{\partial s} + \tilde{J}(s,t,u) \frac{\partial
u}{\partial t} = \tilde{J}(s,t,u) X_{L(s)}(t,u) \textrm{  for } u
\in C^{\infty}( \R \times S^1, M)
\end{equation}
with asymptotic conditions,
\begin{equation}\label{ascond}
\lim_{s \to \pm \infty} u(s,t) = x_{\pm}
\end{equation}
where $x_{-}$ and $x_{+}$ are 1-periodic orbits for $H_2$ and
$H_1$ respectively. Because of the conditions imposed, the
solutions of (\ref{monh12}, \ref{ascond}), stay in a compact set.
Generically the moduli spaces $\mathcal{M}(x_{-}, x_{+})$ are
manifolds of dimension $\mu_{CZ}(x_{+}) - \mu_{CZ}(x_{-})$. Unlike
the solutions of (\ref{eq13}), the solutions of (\ref{monh12}) are
no longer $\R$-invariant and therefore the 0-dimensional moduli
spaces are no longer empty. Notice that the action
$\mathcal{A}_{L(s)} (u(s,\cdot))$ is increasing along a solution
of (\ref{monh12}). Indeed,
\begin{equation}\label{mhact}
\frac{d}{ds} \mathcal{A}_{L(s)} = \|u_s\|_{g_{\tilde{J}(s)}} ^2 -
\int_{S^1} \frac{\partial L}{\partial s} (s,t,u(s,t)) dt \geq 0
\end{equation}
This allows us to define a map, $m$,  between the chain complexes
$$m: CF^a (H_1,J_1) \to CF^a (H_2, J_2)$$
$$m(x_{+}) = \sum_{\mu_{CZ} (x_{+}) = \mu_{CZ} (x_{-})} \#
\mathcal{M}(x_{-}, x_{+}) \langle x_{-} \rangle.$$ The map $m$
preserves the grading and commutes with the differential. It
descends to a morphism in the homology and is called \textit{the
monotonicity homomorphism}, $m(H_1, H_2)$,
\begin{equation}\label{mon}
m(H_1, H_2) : HF_{*} ^{[a,b)} (H_1, J_1) \to HF_{*} ^{[a,b)} (H_2,
J_2)
\end{equation}
\begin{rema}\label{rem1} Standard arguments as in \cite{FH1, CFH} show that
the monotonicity map, $m(H_1, H_2)$, is independent of the choice
of the monotone homotopy used to define it.
\end{rema}
Further, the monotonicity homomorphism satisfies,
\begin{equation}\label{monh}
m(H_2, H_3) \circ m(H_1, H_2)= m(H_1, H_3) \textrm{ for } H_1 \leq
H_2 \leq H_3
\end{equation}
Now we are ready to define the symplectic homology groups of a
nonempty open set $U \subset \C^n$, as the direct limit of the
Floer homology of regular pairs $(H,J)$:
\begin{equation}\label{sh}
S_{*} ^{[a,b)} (U) = \lim_{\longrightarrow} HF_{*} ^{[a,b)} (H,J)
\end{equation}

In what follows, in this subsection, we will outline some results
and constructions concerning the symplectic homology, which will
be important in the proof of our main result. Given $-\infty < a
\leq b \leq c \leq \infty$, we have an exact sequence of chain
complexes given by inclusions,
$$0\longrightarrow C_{*} ^{[a,b)} (H,J) \longrightarrow C_{*} ^{[a,c)}
(H,J) \longrightarrow C_{*} ^{[b,c)} (H,J) \longrightarrow 0$$ and
this generates an exact triangle $\bigtriangleup_{a,b,c} (H,J)$ in
the homology,
$$ HF_{*} ^{[a,b)} (H,J) \longrightarrow HF_{*} ^{[a,c)}
(H,J) \longrightarrow HF_{*} ^{[b,c)} (H,J) \longrightarrow
HF_{*-1} ^{[a,b)} (H,J).$$ $\bigtriangleup_{a,b,c} (H,J)$ commutes
with the monotonicity homomorphism, (\ref{mon}), and gives rise to
an exact triangle, $\bigtriangleup_{a,b,c} (U)$ in symplectic
homology,
$$S_{*} ^{[a,b)} (U) \longrightarrow S_{*} ^{[a,c)}
(U) \longrightarrow S_{*} ^{[b,c)} (U) \longrightarrow S_{*-1}
^{[a,b)} (U).$$ Given triplets $-\infty < a \leq b \leq c \leq
\infty$ and $-\infty < a' \leq b' \leq c' \leq \infty$ with $a
\leq a', b \leq b', c\leq c'$ we consider first the natural map,
given by inclusions,
$$C_{*} ^{[a,b)} (H,J) \longrightarrow C_{*} ^{[a',b')}
(H,J),$$ which gives rise to a map $\sigma$ in homology,
\begin{equation}\label{sigma}
\sigma : HF_{*} ^{[a,b)} (H,J) \longrightarrow HF_{*} ^{[a',b')}
(H,J)
\end{equation}
The map $\sigma$ is compatible with the monotonicity homomorphism
and generates a map $\hat{\sigma}$ in the symplectic homology,
$$\hat{\sigma} : S_{*} ^{[a,b)} (U) \longrightarrow S_{*} ^{[a',b')}
(U).$$ The map $\hat{\sigma}$ commutes with with the triangle
$\bigtriangleup_{a,b,c} (U)$ and generates homomorphisms,
$$\bigtriangleup_{a,b,c} (U) \longrightarrow \bigtriangleup_{a',b',c'}
(U).$$ Given two open and bounded subsets of $\C^n$, $U \subset
V$, we have $\mathcal{H}_{ad}(V) \subset \mathcal{H}_{ad}(U)$.
This observation together with the monotonicity homomorphisms
gives an inclusion morphism, $i_{U,V}$,
\begin{equation}\label{inclus}
i_{U,V} :S_{*} ^{[a,b)} (V) \longrightarrow S_{*} ^{[a,b)} (U)
\end{equation}
For $U \subset V \subset W$, we have,
$$i_{U,W} = i_{U,V} \circ i_{V,W}$$

Consider a regular pair $(H,J)$, and let $c \geq 0$ be a constant.
From (\ref{sigma}) we get a map
$$\sigma(H, c) : HF_{*} ^{[a-c,b-c)} (H,J) \longrightarrow HF_{*} ^{[a,b)}
(H,J)$$ Now observe that the action functionals associated with
$H$ and $H-c$ are related via $\mathcal{A}_{H-c} =
\mathcal{A}_{H}+ c$. This equality translates into an isomorphism,
\begin{equation}\label{isom}
\phi(H-c,H) : HF_{*} ^{[a,b)} (H-c,J) \longrightarrow HF_{*}
^{[a-c,b-c)} (H,J)
\end{equation}
Composing the last two maps, we get a map,
\begin{equation}\label{mon1}
\hat{m}(H-c,H) = \sigma(H,c) \circ \phi(H-c,H) : HF_{*} ^{[a,b)}
(H-c,J) \longrightarrow HF_{*} ^{[a,b)} (H,J)
\end{equation}
On the other hand we have from (\ref{mon}), the monotonicity
homomorphism $m (H-c,H)$. The following lemma, proven in
\cite{He1}, will be useful.
\begin{lm}\label{Hermann}
For any constant $c \geq 0$, $\hat{m}(H-c,H) = m (H-c,H)$.
\end{lm}
We conclude this subsection by outlining a way to compute the
symplectic homology groups for given open set $U$. For this we
need the notion of a cofinal (exhausting) family.
\begin{df} A family of functions $\{H_{\lambda} \}_{\lambda \in
\Lambda}$, where $\Lambda \subset \R$ is unbounded from above, is
called a cofinal family for $U$ if for every $K \in
\mathcal{H}_{ad}(U)$ there exists a number $\lambda'$ s. t.
$H_{\lambda} \geq K$ for $\lambda
> \lambda'$.
\end{df}
Once we have a cofinal family $\{H_{\lambda} \}_{\lambda \in
\Lambda}$, we pair each $H_{\lambda}$ with a compatible almost
complex structure $J_\lambda$. Then one perturbs the family
$(H_{\lambda}, J_{\lambda})$ to get a regular cofinal family or
argues as in \cite{BPS}, Section 4, and the symplectic homology
groups are computed, as,
$$S_{*} ^{[a,b)} (U) = \lim_{\lambda \to \infty} HF_{*} ^{[a,b)}
(H_{\lambda},J_{\lambda})$$ For examples of such computations we
refer to \cite{FHW,CFHW,BPS,CGK,He}.
%%%%%%%%%%%%%%%%%%%%%%%%%%%%%%%%%%%%%%%%%%%%%%%%%%%%%%%%%%%%%%%%%%%%%%%%%%%%%%%%%%%%%%%%%%

\subsection{The definition of the capacity $c_{FH}$.}\label{FlHo}
First recall the definition of a symplectic capacity on $(\R^{2n}
= \C^n, \omega_0=  - d \lambda_0)$.
\begin{df} A symplectic capacity is a map which associates to a
given set $U \subset \C^n$ a number $c(U)$  with the following
properties,
\begin{enumerate}
\item \emph{Monotonicity:} If $U \subset V$ then $c(U) \leq c(V)$,
\item \emph{Symplectic invariance:} $c(\phi(U)) = c(U)$, for any
sympectomorphism $\phi$ of $\C^n$, \item \emph{Homogeneity:} $c(
aU) = a^2 c(U)$ for any real number $a$. \item
\emph{Normalization:} $c (B^{2n} (1)) = c (Z(1)) = \pi$, where
$B^{2n} (1)$ is the unit ball in $\C^n$, centered at the origin
and $Z(1) = \{ z=(z_1, \ldots , z_n) \in \C^n \mid |z_1| < 1 \}$
\end{enumerate}
\end{df}
\begin{rema} Notice that it is sufficient to find such map $c$
with the above properties on open and bounded subsets of $\C^n$,
afterwards we can extend it to any open set as follows,
$$c(U) = \sup \{ c(V) \mid V  \textmd{ is bounded and connected and }
V \subset U \}$$ and to any subset by:
$$c(E) = \inf \{ c(U) \mid U \textmd{ is open and } E \subset U
\}$$
\end{rema}

Now we are ready to review the definition of the Floer-Hofer
capacity as in \cite{He1}. It is based on the computations of the
symplectic homology groups for open balls in \cite{FHW}. We have
\begin{lm} \label{shlm} The symplectic homology groups of an open ball of radius
$R$, $B_R = B^{2n}(R) \subset \C^n$, satisfy
$$S^{[a,b)} _n (B^{2n} (R)) = \Z _2 \textmd{ for }   a
\leq  0 < b \leq \pi R^2, \textmd{ and } 0 \textmd{ otherwise}.$$
$$S^{[a,b)} _{n+1} (B^{2n}(R)) = \Z _2 \textmd{ for } 0 <a \leq \pi R^2 <
b, \textmd{ and } 0 \textmd{ otherwise}.$$
$$S^{[a,b)} _{k} (B^{2n}(R)) = 0 \textmd{ for } k<n \textmd{ or }
n<k <3n$$
\end{lm}
Let $U$ be an open and bounded subset of $\C^n$ and let $r > 0$ be
a number such that $ B^{2n}(r) \subset U$. Pick numbers
$\varepsilon >0 $ such that $\varepsilon < \pi r^2$ and a number
$b > \pi r^2$. Originally, in \cite{FHW}, the following capacity
function was defined. With the inclusion morphism,
$$\sigma _U ^b : S^{[\varepsilon, b)} _{n+1} (U) \to
S^{[\varepsilon,b)}_{n+1} (B^{2n}(R)) = \Z _2$$ we define a
capacity function $c'(U)$ as
$$c'(U) = \inf \{b \mid \sigma _U ^b \textmd{ is onto } \}$$
D. Hermann, was able to extract another capacity from the
symplectic homology which he called the Floer-Hofer capacity and
we adopted his terminology, (see \cite{He1}). Observe that for
large $b$, the natural map,
$$\Z _2 = S^{[0, \varepsilon)} _n ( B^{2n}(\rho)) \to S^{[0, b)}
_n (B^{2n} (\rho))$$ vanishes, (see \cite{V}). Let $R$ be
sufficiently large so that $B_r = B^{2n}(r) \subset U \subset
B^{2n}(R)= B_R$, then we have
$$\begin{CD}
\Z_2 = S^{[0, \varepsilon)} _n (B_R) @> i_R >> SH^{[0,
\varepsilon)} _n(U)@> i_r >> S^{[0, \varepsilon)} _n (B_r) = \Z_2
\end{CD}$$
Since the composition $i_R \circ i_r$ is an isomorphism, it
follows that $0 \neq \alpha_U = i_R(1) \in S^{[0, \varepsilon)} _n
(U)$. One then considers the natural map
$$i^b _U : S^{[0, \varepsilon)} _n (U) \to S^{[0, b)} _n (U)$$
and the Floer-Hofer capacity is defined as
\begin{equation}\label{fh}
c_{FH} (U) = \inf \{ b \mid i^b _U (\alpha_U ) = 0 \}
\end{equation}
The next proposition, relates the capacities $c'$ and $c_{FH}$. It
is proven in \cite{He1} but we sketch a part of the proof for
convenience and better understanding of the nature of the two
capacities.
\begin{prop} The maps $c'$ and $c_{FH}$ are symplectic capacities
and $c' \leq c_{FH}$.
\end{prop}
\textbf{Proof.} Consider the following diagram,
$$\begin{CD}
S_{n+1}^{[0,b)}(B_R) @ >>>S_{n+1}^{[\varepsilon,b)}(B_R) @ >>>
S_{n}^{[0,\varepsilon)}(B_R)=\mathbb{Z}_2 @ >>> S_{n}^{[0,b)}(B_R)\\
@VVV @VVV @VV{i_R}V @VVV\\
S_{n+1}^{[0,b)}(U) @ >>>S_{n+1}^{[\varepsilon,b)}(U) @
>{\partial_U}>>
S_{n}^{[0,\varepsilon)}(U)@ >{i_U ^b}>> S_{n}^{[0,b)}(U)\\
@VVV @VV{\sigma_U ^b}V  @VV{i_r}V @VVV \\
S_{n+1}^{[0,b)}(B_r)=0 @
>>>S_{n+1}^{[\varepsilon,b)}(B_r)=\mathbb{Z}_2 @
 >{\partial_r}>> S_{n}^{[0,\varepsilon)}(B_r)=\mathbb{Z}_2 @ >>>
S_{n}^{[0,b)}(B_r)\\
\end{CD}$$ Here the horizontal arrows are the exact triangles
$\bigtriangleup _{0,\varepsilon, b}$ and the vertical ones are the
inclusion morphisms. We have that $i_r (\alpha_U) = 1$ and
$\partial _r$ is an isomorphism. If $i^b _U (\alpha _U) =0 $ then
there is $\beta \in S_{n+1}^{[\varepsilon,b)}(U)$, such that
$\alpha_U =
\partial_U (\beta)$. We deduce that $\partial_r (\sigma ^b _U
(\beta)) = 1$ and therefore $\sigma ^b _U$ is onto, implying
$c'(U) \leq c_{FH}(U)$. For the fact that $c'$ and $c_{FH}$ are
symplectic capacities we refer to \cite{He1, FHW}.
\begin{rema} D. Hermann, \cite{He}, proves also that the two
capacities are equal on open sets with restricted contact type
boundary. \end{rema}

%%%%%%%%%%%%%%%%%%%%%%%%%%%%%%%%%%%%%%%%%%%%%%%%%%%%%%%%%%%%%%%%%%%%%%%%%%%%%%%%%%%%%%%%%%

\section{Proof of Theorem \ref{D}.} \label{mt} %First recall the definition of
%the Floer-Hofer capacity (in the terminology of D. Hermann,
%\cite{He1}). Denote by $\mathcal{H}_{ad}(U)$, where $U \subset
%\C^n$ is open, the set of all smooth functions $H \in C^{\infty}(
%S^1 \times \C^n)$, which are negative on $\bar{U}$ and $H(t,x) =
%\mu |x|^2$ for $|x|
%>R$ for some $R$ and $\mu \notin \pi \mathbb{Z}$.
%\begin{df} A family of functions $H_{\lambda}$ is called a cofinal
%family for $U$ if for every $K \in \mathcal{H}_{ad}(U)$ there
%exists a number $\lambda_0$ s. t. $H_{\lambda} \geq K$ for
%$\lambda > \lambda_0$.
%\end{df}

Let $\epsilon >0$ be the number given by Lemma \ref{lm1}, we may
assume in addition that $1> \epsilon > 0$. Fix $\epsilon'$ such
that $\epsilon > \epsilon' >0$. For $0 <\tau \leq \epsilon$,
denote by $$V_{\tau} = \psi (N \times B_{\tau} ^k) = \{ x \in U |
\sum_{j=1} ^k z_j ^2 (x) < \tau^2 \}$$ Consider the 1-forms $B_j$
defined on $\R^{2n}$ by $B_j = f \beta_j$, where $f$ is a smooth
function on $\R^{2n}$ such that $f=1$ on $V_{\epsilon'}$ and $f =
0$ on $\R^{2n} \setminus V_{\epsilon}$ and $\beta_j$ are given by
(\ref{beta}). This way we get $k$ one-forms defined on $\R^{2n}$
such that
\begin{equation} \label{B1}
B_j = \beta_j \textmd{ on } V_{\epsilon'}
\end{equation}
and
\begin{equation} \label{B2}
B_j = 0 \textmd{ on } \R^{2n} \setminus V_{\epsilon}
\end{equation}

Now, fix $0 < \delta < \epsilon'$, and consider the set $V_{\delta
/2}$. Using the properties of the capacity $c_{FH}$ we have that
\begin{equation} \label{tau}
c_{FH}(N) \leq c_{FH} (V_{\delta /2}).
\end{equation}

%Now pick a function $K_{\delta}(t,x)$ which generates $\phi$ and
%such that $\parallel K_{\delta}
%\parallel < c_{FH}(N) + \tau(\delta)/8$.

Next we want to construct a cofinal family $H_{\lambda}$ for
$V_{\delta/2}$ for fixed $\delta$. In what follows the parameter
$\lambda$ should be thought as a sufficiently large number since
we will be interested in taking the limit as $\lambda \to \infty$
and so we assume that $\lambda > 16/ \delta$. We mention that the
family we will construct is the one considered by D. Hermann in
\cite{He1}, but adapted for our purposes. Consider smooth
functions $g$ and $h$ on $\R ^{+}$ so that.
\begin{itemize}
\item $h'(t) = \lambda$ for $t \in [\delta /2 + \lambda ^{-1},
\delta /2 + \lambda^{-1/2}]$,

\item $h(t) = - \lambda ^{-1}$ for $t \in [0, \delta /2 - \lambda
^{-1}]$,

\item $h(t) = - \lambda ^{-1} + \lambda^{1/2}$ for $t \geq
\delta/2 + \lambda ^{-1} + \lambda ^{-1/2}$,

\item $h$ is convex on $[\delta/2 - \lambda ^{-1}, \delta/2 +
\lambda^{-1}]$ and concave on $[\delta/2 + \lambda ^{-1/2},
\delta/2 + \lambda ^{-1/2} + \lambda^{-1}]$,

\item $h (\delta/2) <0$

\item $g(t) = - \lambda ^{-1} + \lambda^{1/2}$ for $ t <
(\lambda^{1/6} +1)^2 - \lambda^{-1}$,

\item $g'(t) = \mu/2$ for $t> (\lambda^{1/6} +1)^2$,

\item $g$ is convex on $ [(\lambda^{1/6} +1)^2 - \lambda^{-1},
(\lambda^{1/6} +1)^2]$.

\end{itemize}

Here $\mu \sim \lambda ^{1/6}$ and $\mu \notin \pi \Z$. Now define
$H_{\lambda}$ as follows.
\begin{itemize}

\item $H_{\lambda} (x) = h( \sum _{j=1}^k z_j ^2(x))$ for $x \in
\bigcup_{|\nu |^2< \delta/2 + \lambda ^{-1} + \lambda ^{-1/2}}
N(\nu)$,

\item $H_{\lambda} (x) = g (|x|^2)$ for $ |x| > \lambda ^{1/6}$,

\item $H_{\lambda}(x) = - \lambda ^{-1} + \lambda^{1/2}$ for $x
\in B^{2n} (\lambda ^{1/6}) \setminus \bigcup_{|\nu |^2< \delta/2
+ \lambda ^{-1} + \lambda ^{-1/2}} N(\nu)$.
\end{itemize}
Obviously $H_{\lambda}$ is a cofinal family for $V_{\delta /2}$.
Before we proceed we would like to perturb each $H_{\lambda}$
where it is negative to create non-degenerate critical points. We
do this as follows.  Let $z_0 \in N$, we will create a small
``dimple" at $z_0$. Let $\rho >0$ be such that $B_{\rho}(z_0)
\subset V_{\delta/4}$. Consider a smooth cutoff function $\chi$,
such that $\chi(0) =0$, $\chi(s) =1/2$, for $s \geq \rho /2 $ and
$\chi'(s) >0$ for $s>0$. Denote by $p(x)$ the function
$\chi(r^2(x,z_0))$, where $r(x,z_0)$ is the distance function.
Glue smoothly to $p$ a smooth function $q(x)$ so that $q(x) = 0$
on $B_{\rho/2}(z_0)$ and $q(x) = q_{\lambda}(\sum_{j=1}^k z_j
^2(x))$ for $x$ outside of $V_{\delta/4}$ and $q_{\lambda}$ is a
smooth function on $[\delta/4, \infty)$, such that it is equal to
$1$ on $(\delta/2 - \lambda^{-1}, \infty)$ and $q'_{\lambda} > 0$
on $(\delta/4, \delta/2 - \lambda^{-1})$. Call the new function
$\tilde{q}_{\lambda}$. We assume that it has the following
properties:
\begin{itemize}
\item $1 \geq \tilde{q}_{\lambda} \geq 0$; \item
$\tilde{q}_{\lambda}$ is a Morse function that has global minimum
at $z_0$ equal to $0$; \item the critical points of
$\tilde{q}_{\lambda}$ are contained in $V_{\delta/4}$ (notice that
the gradients $\{\nabla z_j\}_{j=1}^k$ are linearly independent).
\end{itemize}
Now perturb each $H_{\lambda}$ by adding $\lambda^{-2}
(\tilde{q}_{\lambda} (x)-1)$. This way we get a family
$\tilde{H}_{\lambda}(x) = H_{\lambda}(x) + \lambda^{-2}
(\tilde{q}_{\lambda} (x)-1)$. We will abuse the notation and call
the new family $H_{\lambda}$. Again it is a cofinal family for
$V_{\delta /2}$. This way we ensure that for sufficiently small
$\varepsilon >0$ there is a large $\lambda$ so that the only
critical points of $\mathcal{A}_{H_{\lambda}}$ with action in the
interval $[0, \varepsilon)$ are the critical points of
$\tilde{q}_{\lambda} (x)$ which are non-degenerate. We have that
the Conley-Zehnder indices of these critical points, as critical
points of $\mathcal{A}_{H_{\lambda}}$, satisfy:
\begin{equation}\label{CZ}
\mu_{CZ}(x) = m(x)-n \textmd{ for } x \in
Crit(\tilde{q}_{\lambda}) \subset Crit(H_{\lambda})
\end{equation}
where $m(x)$ is the Morse index of $x$, we refer to \cite{SZ}, for
this and other facts concerning the properties of the
Conley-Zehnder index. That is to say that for sufficiently small
$\varepsilon >0$ and large $\lambda$, and $x_0$ - a critical point
of $H_{\lambda}$ with Morse index $l$, then $\Z_2 \langle x_0
\rangle \subset CF_{l-n} ^{[0, \varepsilon)} (H_{\lambda})$. In
particular if $x_0$ is a local minimum of $H_{\lambda}$, then
$\Z_2 \langle x_0 \rangle \subset CF_{n} ^{[0, \varepsilon)}
(H_{\lambda})$.

Next we pair each $H_{\lambda}$ with a compatible almost complex
structure $J_{\lambda}$. We can perturb $J_{\lambda}$ if necessary
to have that the gradient of the function $\tilde{q}_{\lambda}
(x)$ with respect to the metric $g_{J_{\lambda}}$ is Morse-Smale,
see \cite{SZ}, Theorem 8.1. Notice that the critical points of
$\mathcal{A}_{H_{\lambda}}$ may not be non-degenerate. In fact
there are degenerate critical points coming from the region on
which $H_{\lambda} = - \lambda ^{-1} + \lambda ^{1/2}$. In this
situation we can argue as in \cite{BPS}, Section 4, that the
groups $HF_{*} ^{[a,b)} (H_{\lambda}, J_{\lambda})$ are
well-defined as long as $a,b \notin \Sigma(H_{\lambda})$, see
especially Remark 4.4.1 in \cite{BPS}.

Consider the function $\hl(t,x)$, defined as follows,
\begin{equation}\label{}
\hl (t,x) = 0 \textmd{ for } 0 \leq t <1/2; \textmd{ and }
\hl(t,x) = 2 H_{\lambda}(x) \textmd{ for } 1/2 \leq t < 1
\end{equation}
If we consider the action functional, associated with $\hl$, it
has the form:
$$\mathcal{A}_{\hl} = - \int_{S^1} x^* \lambda_0 - \int_{1/2} ^1 2
H_{\lambda}(x(t)) dt.$$ Straightforward computations show that
$\mathcal{A}_{H_{\lambda}}$ and $\mathcal{A}_{\hl}$ have the same
critical points with the same critical values and Conley-Zehnder
indices. In fact, they generate the same Floer homology groups.
\begin{prop}\label{equiv} $HF_{*} ^{[a,b)} (H_{\lambda}) \cong HF_{*}
^{[a,b)} (\hl)$ for all $-\infty < a \leq b \leq \infty$.
\end{prop}
Observe that the function $\hl$ is not smooth. Despite that it has
well defined Floer homology. The reason is that the set up for the
Floer homology involves Sobolev spaces of the type $W^{1,p}$ and
all the analysis is carried over initially in a weak sense and
then elliptic ``bootstrapping" arguments are applied for the
smoothness of the solutions. The same type of analysis can be
carried for piecewise smooth functions. Besides, the critical
points of $\mathcal{A}_{\hl}$ are smooth loops. So, in a way Floer
homology ``forgives" slight irregularities of the Hamiltonians.
Now observe that the functions $H_{\lambda}$ and $\hl$ generate
the same time 1-maps. In that situation Proposition \ref{equiv} is
a consequence of the discussion in \cite{Se}, Section 4.

Next pick a compactly supported Hamiltonian function $K_{\delta}$,
which generates $\phi$ and such that $\|K_{\delta} \| < c <
c_{\delta} = c_{FH}(V_{\delta /2})$, where $c$ is some positive
number satisfying the previous inequality, see (\ref{tau}). Denote
by $\hk$ the following function
\begin{alignat*}{4}
H_{\lambda} \natural K_{\delta} (t, x) = 2(K_{\delta}(2t, x) -
\sup K_{\delta}) \textmd{ for }& 0 \leq t <1/2; \\ \hk(t,x) = 2
H_{\lambda}(x) \textmd{ for } 1/2 \leq t < 1.
\end{alignat*}
We will be interested in the critical points of the action
functional associated with $\hk$.
$$\mathcal{A}_{\hk} (x) = - \int_{S^1} x^* \lambda_0 - 2 \int_0 ^{1/2}
(K_{\delta}(2t,  (x)) - \sup K_{\delta}) dt - 2 \int_{1/2}^1
H_{\lambda} (x) dt$$ To be more precise we are going to show that
this functional possesses a critical point (i.e., a 1-periodic
orbit of $\hk$) with action in the interval $[0, c_{\delta}]$, for
sufficiently large $\lambda$. Observe that any critical point of
$\mathcal{A}_{\hk}$ consists of two arcs, one is a trajectory of
the flow of $K_{\delta}$, followed by a trajectory of
$H_{\lambda}$. Notice that we have
$$\hl - c \leq \hk \leq \hl$$
and
$$\mathcal{A}_{\hl} \leq \mathcal{A}_{\hk} \leq
\mathcal{A}_{\hl}+c$$ The next lemma is a modification of
Corollary 5.9 in \cite{He}, but notice that we assume less in our
case.
\begin{lm} Let $c < c_{\delta}$ be as above. Then for sufficiently
large $\lambda$, $\hk$ has a 1-periodic orbit with action in the
interval $[0,c_{\delta}]$.
\end{lm}
\textbf{Proof :} Pick a sufficiently small $\varepsilon$ so that
$0 < \varepsilon < c_{\delta} -c$. Let $B_r$ and $B_R$ be balls
centered at $z_0 \in N \subset \C^n$ with radii $r$ and $R$
respectively so that $B_r \subset V_{\delta/2} \subset B_R$. We
know from Lemma \ref{shlm} that
\begin{equation}
S_n^{[-c, \varepsilon)}(B_r) \simeq \Z_2 \simeq S_n^{[-c,
\varepsilon)}(B_R)
\end{equation}
Moreover we can easily construct cofinal families for $B_r$ and
$B_R$, respectively with a single ``dimple", i.e. unique local
minimum at $z_0$ for both families in the spirit of what we did
with $H_{\lambda}$. Our arguments above show that then the
generator of the symplectic homology groups $S_n^{[-c,
\varepsilon)}(B_r)$ and $S_n^{[-c, \varepsilon)}(B_R)$ is the
class of the constant solution, i.e. $[z_0]$. Consider the
following diagram for sufficiently large $R$,
\begin{equation}\label{d2}
\begin{CD}
\Z_2 = S_{n}^{[-c,\varepsilon)}(B_R) @
>{\sigma_R}>>S_{n}^{[0, c + \varepsilon)}(B_R) = \Z_2 \\
@VVV @VVV \\
\{0\} \neq S_{n}^{[-c,\varepsilon)}(V_{\delta/2}) @
>{\sigma '}>>S_{n}^{[0, c + \varepsilon)}(V_{\delta/2}) \neq \{0\}\\
@VVV  \\
\Z_2 = S_{n}^{[-c ,\varepsilon)}(B_r)\\
\end{CD}
\end{equation}
where the vertical arrows are the inclusion morphisms and the fact
that $$S_{n}^{[0, c + \varepsilon)}(V_{\delta/2}) \neq \{0\},$$
follows from the definition of the Floer-Hofer capacity,
(\ref{fh}). This diagram, (together with the definition of the
Floer- Hofer capacity), implies that the map $\sigma '$ is
nonzero. In fact, the map $\sigma '$ keeps ``alive" the class of
$z_0$.  Next consider the commutative diagram.
\begin{equation}\label{d3}
\begin{CD}
HF_{n}^{[-c,\varepsilon)}(H_{\lambda}, J_{\lambda}) @
>{\sigma(H_{\lambda}, c)}>>HF_{n}^{[0, c + \varepsilon)}(H_{\lambda}, J_{\lambda}) \\
@VVV @VVV \\
\{0\} \neq S_{n}^{[-c,\varepsilon)}(V_{\delta/2}) @
>{\sigma '}>>S_{n}^{[0, c + \varepsilon)}(V_{\delta/2}) \neq \{0\}\\
\end{CD}
\end{equation}
Here the vertical arrows are the direct limit morphisms, which are
surjective for sufficiently large $\lambda$. This diagram implies
that the map ${\sigma(H_{\lambda}, c)}$ must be nonzero. But then
from  (\ref{mon1}), we have that the map
$$\hat{m}(H_{\lambda} -c, H_{\lambda}) : HF_n ^{[0,c+\varepsilon)}
(H_{\lambda} -c, J_{\lambda}) \rightarrow HF_n
^{[0,c+\varepsilon)} (H_{\lambda}, J_{\lambda})$$ is nonzero. From
Lemma \ref{Hermann}, we know that $\hat{m}(H_{\lambda} -c,
H_{\lambda})= m(H_{\lambda} -c, H_{\lambda})$. And this shows that
$m(H_{\lambda} -c, H_{\lambda}) \neq 0$. Denote by $\bar{m}(\hl
-c, \hl)$ the monotonicity map between $HF_n ^{[0,c+\varepsilon)}
(\hl -c, J_{\lambda})$ and $HF_n ^{[0,c+\varepsilon)} (\hl,
J_{\lambda})$. It is not hard to see that $\bar{m}$, agrees with
the map induced by $m$ through the isomorphism of Proposition
\ref{equiv}. Then this map is nonzero. If $\hk$ did not have a
1-periodic orbit with action in $[0,c_{\delta}]$, then the Floer
homology group $HF_n ^{[0,c + \varepsilon)} (\hk, J_{\lambda})$
would be well-defined and equal to zero. But then the monotonicity
map $\bar{m}(\hl -c, \hl)= \bar{m}(\hl -c, \hk) \circ \bar{m}(\hk,
\hl)$, would have been zero, which is a contradiction. $\Box$

The previous Lemma gives the existence of a 1-periodic orbit,
$x_{\lambda} (t)$ of $\hk$ with bounded action for sufficiently
large $\lambda$ and moreover that this solution is a deformation
of the constant class of $0 \neq [z_0] \in HF_n ^{[0,c+
\varepsilon)} (\hl -c, J_{\lambda})$ or in other words,
\begin{equation}
0 \neq \bar{m}(\hl -c, \hk) ([z_0]) \in HF_n ^{[0, c +
\varepsilon)} (\hk, J_{\lambda})
\end{equation}
and
\begin{equation}\label{mono}
x_{\lambda} = \bar{m}(\hl -c, \hk) (z_0)
\end{equation}
This observation will be important later on. The periodic orbit
$x_{\lambda} (t)$, satisfies the equations,
\begin{equation}\label{arcs}
\begin{split}
\dot{x}_{\lambda}(t) = 2 X_{K_{\delta}(2t)} (x(t)) \textmd{ for }
t
\in (0, 1/2)\\
\dot{x}_{\lambda}(t) = 2 X_{H_{\lambda}} (x(t)) \textmd{ for } t
\in ( 1/2,1)
\end{split}
\end{equation}
Denote by $\phi_{\delta}^t$ the flow of $K_{\delta}$ and by
$\varphi_{\lambda}^t$ the flow of $H_{\lambda}$.

Next we claim that for sufficiently large $\lambda$, $x_{\lambda}
(0) \in N_{\nu}$, where $\nu = (\nu_1, \ldots \nu_k)$ and
$\sum_{j=1} ^k \nu_j ^2 < \delta /2 + \lambda ^{-1} +\lambda
^{-1/2}$. Indeed, if we assume that this is not the case then we
have two possibilities: either $ x_{\lambda} (t) = x_{\lambda}
(1/2)$ for $t \in [1/2,1]$ or (perhaps after taking sufficiently
large $\lambda$ so large that the ball $B(\lambda ^{1/6}) \supset
\textmd{supp} K_{\delta}$), $ x_{\lambda}(t) = x_{\lambda}(0)$ for
$t \in [0,1/2]$. In the former case we have that $x_{\lambda}(t) =
\phi_{\delta} ^t (x_{\lambda}(0))$ is a 1-periodic solution for
$K_{\delta}$, then its action satisfies,
\begin{equation} \label{action1}
\mathcal{A}_{\hk} (x_{\lambda}(t)) =
\mathcal{A}_{K_{\delta}}(x_{\lambda}(t)) - \sup K_{\delta} -
C(\lambda)
\end{equation}
where $C(\lambda) = \lambda^{1/2} - \lambda ^{-1}$. Since
$K_{\delta}$ is compactly supported, the critical values of
$\mathcal{A}_{K_{\delta}}(x(t)) - \sup K_{\delta}$ are bounded and
therefore for large $\lambda$, the right hand side of
(\ref{action1}) will be very negative, which is a contradiction
with the fact that $\mathcal{A}_{\hk}(x_{\lambda}(t)) \geq 0$.
Similarly in the latter case, $x_{\lambda} (t)$ is a 1-periodic
orbit for $H_{\lambda}$, satisfying the equation $-i
\dot{x}_{\lambda} = \mu x_{\lambda}$. Then we have for the action
of the periodic orbit $x_{\lambda}$,
\begin{equation}\label{action2}
\mathcal{A}_{\hk} (x_{\lambda}) \leq \int_0 ^1 (\mu |x_{\lambda}
(t)|^2 - g(|x_{\lambda}|^2))dt + c \leq \mu (\lambda ^{1/6} +1 )^2
- C(\lambda) +c
\end{equation}
where $C(\lambda)$ is given before. Because of our choice of $\mu$
and since $c$ is bounded, for large $\lambda$ this action will be
very negative, which is a contradiction.

Now fix a very large $\lambda = \lambda(\delta)$ so that
$B(\lambda ^{1/6}) \supset V_{\epsilon}$, $\lambda ^{-1} + \lambda
^{-1/2} < \delta /2$ and the 1-periodic orbit of $\hk$,
$x_{\lambda}$, satisfies $x_{\lambda}(0) \in N(\nu)$ for some $\nu
= \nu (\delta) = (\nu_1, \ldots, \nu_k)$, with $(\sum_{j=1} ^k
\nu_j ^2)^{1/2} < \delta/2 + \lambda ^{-1} + \lambda ^{-1/2} <
\delta$. From (\ref{arcs}), follows that we can write
$x_{\lambda(\delta)} (t) = \phi_{\delta}^{2t} (x_{\lambda(\delta)}
(0))$, for $t \in [0,1/2]$ and $x_{\lambda(\delta)} (t)=
\varphi_{\lambda}^{2t-1}(x_{\lambda(\delta)} (1/2))$, for $t \in
[1/2,1]$. Then we have,
$$\phi_{\delta}^1 (x_{\lambda(\delta)} (0)) = x_{\lambda(\delta)}
(1/2) = (\varphi_{\lambda}^1)^{-1}(x_{\lambda(\delta)} (0))$$ We
argue that $(\varphi_{\lambda}^1)^{-1}(x_{\lambda(\delta)} (0))
\in \mathcal{L}_{x_{\lambda(\delta)} (0)} N(\nu(\delta))$. Indeed,
if $x \in N(\nu(\delta))$ then the flow $\varphi_{\lambda}^t x$ of
    $H_{\lambda}$, satisfies an equation of the form (\ref{leaf}),
with coefficients $$\lambda_j = 2 h'(\sum_{i=1}^k \nu_i ^2) \nu_j
= 2\lambda(\delta) \nu_j$$ and therefore the flow of $H_{\lambda}$
is on the leaf through $x$. The flow $(\varphi_{\lambda}^t)^{-1}$
is generated by the Hamiltonian $\tilde{H}_{\lambda} (x) = -
H_{\lambda}(\varphi_{\lambda}^t x)$. From this it is not hard to
see that the flow $(\varphi_{\lambda}^t)^{-1}$ on $N(\nu(\delta))$
satisfies an equation of the form,
\begin{equation}\label{leaf1}
\dot{x}(t) = \sum_{j=1} ^k \gamma_j X_{z_j} (x(t))
\end{equation}
and this shows that
$(\varphi_{\lambda}^t)^{-1}(x_{\lambda(\delta)} (0)) \in
\mathcal{L}_{x_{\lambda(\delta)} (0)} N(\nu(\delta))$ for any $t$
and in particular for $t=1$. To summarize we demonstrated that
$x_{\lambda(\delta)} (0) \in N(\nu(\delta))$ satisfies
$$\phi(x_{\lambda(\delta)} (0)) =\phi_{\delta}^1
(x_{\lambda(\delta)} (0)) \in \mathcal{L}_{x_{\lambda(\delta)}
(0)} N(\nu(\delta))$$ The next lemma is crucial since it will
allow us to take a limit as $\delta$, respectively
$\lambda(\delta)$ goes to $0$.
\begin{lm}
The length of the arc
$l((\varphi_{\lambda}^t)^{-1}(x_{\lambda(\delta)} (0))|_{t \in
[0,1]})$ is bounded independently of $\delta$.
\end{lm}
\textbf{Proof :} In view of (\ref{leaf1}) this statement is
equivalent to showing that each of the coefficients $\gamma_j$, $j
=1, \ldots, k$ is uniformly bounded. Recall from (\ref{mono}),
that the periodic orbit $x_{\lambda(\delta)} (t)$ is a deformation
of the constant solution of $\hl - c$ through a monotone homotopy.
From Remark \ref{rem1}, we know that the map $\bar{m}(\hl -c,
\hk)$ is independent of the choice of the monotone homotopy of
Hamiltonians used to define it. This allows us to choose a
particular regular monotone homotopy $(L(s), \tilde{J}(s))$ which
realizes $\bar{m}(\hl -c, \hk)$. We pick $L$ of the form
\begin{equation}
L(s, t, u(s,t)) = (1 - \kappa(s)) (\hl(t,u) - c) + \kappa(s) \hk
(t,u)
\end{equation}
where $\kappa(s)$ is a smooth function on $\R$ so that $\kappa(s)
= 0$ for $s \leq - s_0$; $\kappa(s) = 1$, for $s \geq s_0$ and
$\kappa '(s) \geq 0$ on $(-s_0, s_0)$. Of course we assume that
$\tilde{J} (s,t)$ is a regular homotopy of families of almost
complex structures so that $\tilde{J} (s,t) =
J_{\lambda(\delta)}(t)$, for $s \in (-\infty, -s_0] \bigcup [s_0,
\infty)$. Consider now the equation,
$$u_s + \tilde{J}(s,t,u(s,t)) ( u_t - X_{L(s,t,u(s,t)}) =0$$
Our arguments imply that it possesses a solution $u(s,t)$, such
that
$$\lim_{s \to -\infty} u(s,t) = z_0$$
and
$$\lim_{s \to \infty} u(s,t) = x_{\lambda(\delta)} (t)$$
In view of this and integrating (\ref{mhact}) over $\R$, we get,
for our particular case, the following inequality,
\begin{equation}\label{action3}
\mathcal{A}_{H_{\lambda} - c}(z_0) -
\mathcal{A}_{\hk}(x_{\lambda(\delta)} (t))\geq \frac{1}{2}
\int_{\R \times S^1} (\|u_s\|_{g_{\tJ (s)}}^2 + \|u_t -
X_{L}\|_{g_{\tJ (s)}}^2)dsdt
\end{equation}
The left-hand side of (\ref{action3}) is bounded from above by
$c_{\delta} = c_{FH}(V_{\delta/2})$. We are going to work with the
right-hand side. Recall that $X_{\hk} =2 X_{H_{\lambda}}$ for $t
\in (1/2 , 1)$. Using this we get the following inequality for the
right-hand of (\ref{action3}),
\begin{equation}\label{ac1}
\begin{split}
\frac{1}{2} \int_{\R \times S^1} (\|u_s\|_{g_{\tJ (s)}}^2 + \|u_t
- X_{L}\|_{g_{\tJ (s)}}^2)dsdt\\ \geq  \frac{1}{2} \int_{\R \times
[1/2,1]} %\int_{1/2} ^1
 (\|u_s\|_{g_{\tJ (s)}}^2 + \|u_t - X_{L}\|_{g_{\tJ (s)}}^2)dsdt
\\ = \frac{1}{2} \int_{\R} \int_{1/2}^1 (\|u_s\|_{g_{\tJ (s)}}^2 +
 \|u_t - 2X_{H_{\lambda}}\|_{g_{\tJ (s)}}^2)dsdt
\end{split}
\end{equation}
Now recall the 1-forms $B_j$, $j=1, \ldots, k$, which we
introduced in the beginning of this section, see
(\ref{B1},\ref{B2}). We claim that $B_j (X_{H_{\lambda}}, \cdot )
= 0$. This is easy to be seen, to be the case on $\R^{2n}
\setminus V_{\epsilon}$, since there $B_j = 0$. On $V_{\epsilon}
\setminus V_{\epsilon'}$, it is true because $H_{\lambda} =
\textmd{const}$, there. On $V_{\epsilon'}$, $B_j = \beta_j$, and
on that region $X_{H_{\lambda}}$ is a linear combination of
$\{X_{z_j}\}_{j=1} ^k$, and our claim follows from (\ref{bet}).
Choose a constant $C_1 >0$ so that for $j=1, \ldots , k$ and all
$\xi, \eta \in \R^{2n}$ we have
\begin{equation}\label{C1}
|dB_j (\xi, \eta)| \leq C_1 |\xi| |\eta|
\end{equation}
Consider the space of all almost complex structures $J$ on $\C^n$,
compatible with $\omega_0$. Denote, as before, by $g_J$ the
corresponding metric, i.e., $g_J(\cdot , \cdot) = \omega_0 (\cdot,
J \cdot)$. Since the set $\bar{V}_\epsilon$ is a compact subset of
$\C^n$, there is a constant $C_2 >0$ so that on $\bar{V}_\epsilon$
we have that,
$$\|\xi \|_{g_J} \geq \sqrt{C_2} \|\xi \|_{g_{J_0}} = \sqrt{C_2}
|\xi|$$ for any $\xi \in \C^n$. Here $J_0 = i$ is the standard
complex structure on $\C^n$. In view of our discussion above and
(\ref{action3}, \ref{ac1}) we obtain,
\begin{align*}
c_{\delta}  &\geq  \mathcal{A}_{H_{\lambda} - c}(z_0) -
\mathcal{A}_{\hk}(x_{\lambda(\delta)} (t))\\ &\geq \frac{1}{2}
\int_{\R \times S^1} (\|u_s\|_{g_{\tJ (s)}}^2 + \|u_t -
X_{L}\|_{g_{\tJ (s)}}^2)dsdt \\ &\geq  \frac{1}{2} \int_{\R \times
[1/2,1]} (\|u_s\|_{g_{\tJ (s)}}^2 + \|u_t - 2
X_{H_{\lambda}}\|_{g_{\tJ (s)}}^2)dsdt \\ &\geq  \frac{1}{2}
\int_{\R \times [1/2,1]}^{*} (\|u_s\|_{g_{\tJ (s)}}^2 + \|u_t - 2
X_{H_{\lambda}}\|_{g_{\tJ (s)}}^2)dsdt \\ &\geq  \int_{\R \times
[1/2,1]}^{*} (\|u_s\|_{g_{\tJ (s)}} \|u_t - 2
X_{H_{\lambda}}\|_{g_{\tJ (s)}})dsdt \\ &\geq  \frac{C_2}{C_1}
\int_{\R \times [1/2,1]}^{*} |dB_j (u_s , u_t - 2
X_{H_{\lambda}})| ds dt \\ &=  \frac{C_2}{C_1} \int_{\R \times
[1/2,1]}^{*} |dB_j (u_s , u_t )|ds dt = \frac{C_2}{C_1} \int_{\R
\times [1/2,1]} |dB_j (u_s , u_t )|ds dt \\ &\geq   \frac{C_2}{C_1}
| \int_{\R \times [1/2,1]} dB_j (u_s , u_t )ds dt| =
\frac{C_2}{C_1}|\int_{1/2} ^1  x_{\lambda(\delta)} (t)^* B_j  dt -
\int_{1/2} ^1 z_0 ^* B_j  dt| \\&= \frac{C_2}{C_1} |\int_{1/2} ^1
x_{\lambda(\delta)} (t) ^* B_j dt| = \frac{C_2}{C_1}|\gamma_j|
\end{align*} In the above formulas the last couple of equalities follow from
Stokes' Theorem and (\ref{leaf1}) and $\int ^*$ means integrating
over the part of the trajectory which is contained in
$\bar{V}_\epsilon$. So far, we obtained that for each $j= 1,
\ldots, k$, the coefficients $\gamma_j$ are bounded by $c_{\delta}
C_1 / C_2$. Notice that $c_{\delta} \leq c_{FH} (V_{\epsilon})$
and so it is bounded by a constant independent of $\delta$ and so
are the coefficients $\gamma_j$. All this shows that the length of
the arc $l(x_{\lambda(\delta)}(t)|_{t \in [1/2,1]})$ is bounded
independently of $\delta$. $\Box$

Repeating the arguments above for any $\delta \in (0, \epsilon')$
and applying the Arzela-Ascolli Theorem, we can find a sequence
$\{\delta_m\}_{m=1} ^{\infty}$ converging to $0$ so that
$$\lim_{m \to \infty} x_{\lambda(\delta_m)}(0) = x_0 \in N$$
$$\phi(x_{\lambda(\delta_m)}(0)) \to \phi(x_0)$$
and $x_0$ and $\phi(x_0)$ are connected by an arc which is
contained in the leaf $\mathcal{L}_{x_0} N$. This proves Theorem
\ref{D}. $\Box$

%%%%%%%%%%%%%%%%%%%%%%%%%%%%%%%%%%%%%%%%%%%%%%%%%%%%%%%%%%%%%%%%%%%%%%%%%%%
%%%%%%%%%%%%%%%%%%%%%%%%%%%%%%%%%%%%%%%%%%%%%%%%%%%%%%%%%%%%%%%%%%%%%%%%%%%

\section{Proof of Theorem \ref{DD}.} \label{appl}
Theorem \ref{DD} is a consequence of Theorem \ref{D} and the
following lemmata.
\begin{lm} The level submanifold $N_{(c, c_1, \ldots , c_{k-1})}$
is of $k$-contact type in $(\R^{2n}, \omega_0)$.
\end{lm}
\textbf{Proof :} First we notice that a symplectic change of
coordinates does not change the property of a submanifold to be of
$k$- contact type. Making a symplectic change of coordinates$(x_j,
y_j) \to (I_j, \theta_j)$ , where $I_j = r_j ^2/2$ and as before
$x_j -i y_j = r_j e^{i \theta_j}$. In these coordinates $\omega _0
= d \alpha_0$ with $\alpha_0 = -I_1 d \theta_1 - \ldots -I_n d
\theta_n$ and $N_{(c, c_1, \ldots , c_{k-1})} = \{(I_j, \theta_j)|
I_1 = c_1/2, \ldots I_{k-1} = c_{k-1}/2, \sum_{j=k}^n m_j I_j = c
- 1/2 \sum_{j=1}^{k-1} m_j c_j\}$. For $1 \leq j \leq k-1$,
consider the one-forms $\alpha_j = \alpha_0 - d \theta_j$.
Obviously we have $d \alpha_j = \omega_{0}$ for $0 \leq j \leq
k-1$. Next we see that $Ker \omega_0 |_N = span ( X_{0} =
\sum_{j=1} ^n - m_j \frac{\partial}{\partial \theta_j}, X_{1} =
-\frac{\partial}{\partial \theta_1}, \ldots , X_{k-1}=
-\frac{\partial}{\partial \theta_{k-1}})$. We want to show that
the restrictions of $\alpha_0, \alpha_1, \ldots \alpha_{k-1}$ to
$Ker \omega_0 |_N$ are linearly independent on $N= N_{(c, c_1,
\ldots , c_{k-1})}$ and we check that on the basis $X_{0}, X_{1},
\ldots , X_{k-1}$. Denote by $A$ the $k \times k$ matrix with
entries $a_{i, j} = \alpha_{j-1} (X_{i-1})$ for $1 \leq i, j \leq
k$. Then we have that $a_{1,1} = c$ and $a_{1, j} = c +  m_{j-1}$
for $2 \leq j \leq k$, $a_{i,j} = c_{i-1}/2 + \delta_{i,j}$ for $2
\leq i,j \leq k$ where $\delta_{i,j}$ denotes the Kronecker
symbol. It is not hard to compute that $ \det A = c- \frac{1}{2}
\sum_{j=1}^{k-1} m_j c_j > 0$, because of our assumption
(\ref{cond}). This completes the proof of the lemma.

\begin{lm} The Floer-Hofer capacity, $c_{FH} (N_{r_1, \ldots
r_n}^k) = \min_p \{ \pi r_p ^2\}$, where $N_{r_1, \ldots r_n}^k
=\big\{ \mid z_1 \mid =r_1, \ldots \mid z_{k-1} \mid =r_{k-1},
\sum_{j=k}^n \frac{\mid z_j\mid ^2}{r_j ^2} =1 \big\}$ and $z_j =
x_j + i y_j$.
\end{lm}
\textbf{Proof :} First observe that $c_{FH}(N_{r_1, \ldots r_n}^k)
\leq \min_p \{ \pi r_p ^2\}$. Indeed we have that $N_{r_1, \ldots
r_n}^k) \subset Z_{r_j}$, for $j=1, \ldots n$, where $Z_{r_j} = \{
z \in \C^n | |z_j| < r_j \}$ and the claim follows from the
properties of the capacity.

Next we are going to argue that $c_{FH} (N_{r_1, \ldots r_n}^k)
\geq \min_p \{ \pi r_p ^2\}$. For this we use arguments similar to
those in \cite{FHW}, where the symplectic homology of ellipsoids
and polydisks is computed. Because of that we will be somewhat
sketchy. Essentially the idea is to exploit the product structure
of $(\R^{2n} = \C^n, \omega_0)$. For sufficiently small
$\varepsilon > 0$ consider a neighborhood $V_{\varepsilon}$ of
$N_{r_1, \ldots r_n}^k$ of the form,
\begin{equation*}
V_{\varepsilon} = \bigg\{ z \in \C^n | 1-  \varepsilon <\frac{|z_j|^2}{
r_j ^2} < 1+ \varepsilon \textmd{ for } j=1, \ldots, k-1 \textmd{
and } 1-  \varepsilon < \sum_{j=k}^n \frac{|z_j|^2}{ r_j ^2} < 1+
\varepsilon \bigg\}
\end{equation*}
For $V_{\varepsilon}$ we are going to build a cofinal family of
Hamiltonians of the form:
$$H_{\lambda} (z_1, \ldots, z_n) = \sum_{j=1}^{k-1} \rho_{\lambda} \bigg(\frac{|z_j|^2}{
r_j ^2}\bigg) + \rho_{\lambda} \bigg(\sum_{j=k}^n \frac{|z_j|^2}{ r_j ^2}\bigg)$$
where the functions $\rho_{\lambda}: \R \to \R$ satisfy,
\begin{itemize}
\item $\rho_{\lambda}$ is symmetric with respect to 1, i.e.
$\rho_{\lambda}(1+s) = \rho_{\lambda}(1-s)$ and has unique
non-degenerate minimum at $1$; \item $\rho_{\lambda}'(s) =
\rho_{\lambda}'(\infty) = const$ for $s \geq s_0(\lambda) >>1$;
\item $\rho_{\lambda}'(s) > 0$ for $s>1$; \item $\rho_{\lambda}(s)
<0$ for $s \in [1 -\varepsilon, 1+ \varepsilon]$; \item for each
$\lambda$ the equations $-i \dot{z} = \rho_{\lambda}'(\infty) z$
have no non-trivial 1-periodic solutions; \item for $\lambda >
\lambda'$, $\rho_{\lambda} > \rho_{\lambda'}$.
\end{itemize}
Then one perturbs perturbs $H_{\lambda}$ by small perturbation
$\Delta_{\lambda}$ so that $H_{\lambda}+\Delta_{\lambda} \in
\mathcal{H}_{reg} (V_{\varepsilon})$ and the actions of 1-periodic
orbits of $H_{\lambda}+\Delta_{\lambda}$ are near the actions of
the 1-periodic orbits of $H_{\lambda}$. We abuse the notation and
denote the perturbed family again by $H_{\lambda}$. Then
Proposition 5, in \cite{FHW}, tells us that a minimal non-negative
action of periodic orbit of $H_{\lambda}$, of Conley-Zehnder index
$n+1$ will be greater than $\min_p \{ \pi r_p ^2 (1-\varepsilon)\}
-\tau(\lambda)$ for some $\tau(\lambda) >0$ and such that
$\lim_{\lambda \to \infty} \tau(\lambda) =0$. This immediately
gives us,
$$c_{FH} (V_{\varepsilon}) \geq \min_p \{ \pi r_p ^2
(1-\varepsilon)\}.$$ Passing to the limit as $\varepsilon \to 0$
we get $$c_{FH} (N_{r_1, \ldots r_n}^k) \geq \min_p \{ \pi r_p
^2\}.$$ This completes the proof of the lemma.$\Box$

\textbf{Proof of Theorem \ref{DD}:} Denote by $\phi$ the time-one
map of the Hamiltonian $H_0$ given by (\ref{harmosc}), and by
$\psi$ the time-one map of $H_0 + H_1$. Since $\phi (\mathcal{L}_N
(x)) = \mathcal{L}_N (x)$, we have to show that there exists $x
\in N$ such that
$$\phi ^{-1} \circ \psi (x) \in \mathcal{L}_N (x)$$
The map $\phi ^{-1} \circ \psi$ is the time-one map of the flow
generated by the Hamiltonian $H_0 (\psi^t (x)) + H_1(t, \psi^t
(x)) - H_0 (\psi^t (x)) = H_1(t, \psi^t (x))$. By the preceding
lemmata and the properties of the capacity $c_{FH}$, we know that
$N_{c,c_1, \ldots, c_{k-1}}$ is of $k$-contact type and
$$c_{FH}(N_{c,c_1, \ldots, c_{k-1}})=\min \bigg\{ \min_{p=1,\ldots, k-1}
\{\pi c_p \}, \min_{p=k,\ldots, n} \{\pi
\frac{(2c-\sum_{j=1}^{k-1} m_j c_j)}{m_p}\}\bigg\}$$ Thus we have
$E(\phi ^{-1} \circ \psi) \leq \parallel H_1(t, \psi^t (x))
\parallel < c_{FH}(N_{c,c_1, \ldots, c_{k-1}})$. Now Theorem
\ref{D} yields easily Theorem \ref{DD}. $\Box$

%%%%%%%%%%%%%%%%%%%%%%%%%%%%%%%%%%%%%%%%%%%%%%%%%%%%%%%%%%%%%%%%%%%%%%%%%%%%%%%

%%%%%%%%%%%%%%%%%%%%%%%%%%%%%%%%%%%%%%%%%%%%%%%%%%%%%%%%%%%%%%%%%%%%%%%%%%%
\subsection*{Acknowledgments}I would like to thank P. Albers, T. Ekholm, D. Hermann, K. Honda, H. Hofer and K.
Wysocki for the stimulating discussions and the interest in this
paper. Parts of this work were done during the author's visits at
The University of Melbourne and FIM of ETH Z\"urich and he wishes
to acknowledge the hospitality. Last but not least, the author
expresses his gratefulness for the stimulating scientific
environment to the faculty of the Mathematics Department at the
University of Southern California where he held the position of
Busemann Assistant Professor from August, 2002 until August, 2005.
%%%%%%%%%%%%%%%%%%%%%%%%%%%%%%%%%%%%%%%%%%%%%%%%%%%%%%%%%%%%%%%%%%%%%%%%%%%

%%%%%%%%%%%%%%%%%%%%%%%%%%%%%%%%%%%%%%%%%%%%%%%%%%%%%%%%%%%%%%%%%%%%%%%%%%%
%%%%%%%%%%%%%%%%%%%%%%%%%%%%%%%%%%%%%%%%%%%%%%%%%%%%%%%%%%%%%%%%%%%%%%%%%%%

\end{document}